\documentclass [11pt,oneside]{article}

\reversemarginpar \textwidth=14cm \textheight=19cm

\usepackage{amssymb} \usepackage{amsfonts} \usepackage{amsmath}
\usepackage{amsthm} \usepackage{epsfig} 

\usepackage{caption}

\newcommand{\mettifig}[1]{\epsfig{file=#1}}
\newcommand{\lgh}[1]{{\rm lgh}(#1)}
\newcommand{\Sol}{{\rm Sol}}
\newcommand{\chiorb}{\chi^{\rm orb}}

\newtheorem{lemma}{Lemma}[section]
\newtheorem{teo}[lemma]{Theorem}
\newtheorem{rem}[lemma]{Remark} 
\newtheorem{prop}[lemma]{Proposition}
\newtheorem{cor}[lemma]{Corollary}

\newcommand{\matR} {\ensuremath {\mathbb{R}}}

\newcommand{\matZ} {\ensuremath {\mathbb{Z}}}

\newcommand{\matP} {\ensuremath {\mathbb{P}}}
\newcommand{\matH} {\ensuremath {\mathbb{H}}}
\newcommand{\matS} {\ensuremath {\mathbb{S}}}
\newcommand{\matE} {\ensuremath {\mathbb{E}}}
\newcommand{\matRP} {\ensuremath {\mathbb{RP}}}

\newcommand{\calN} {\ensuremath {\mathcal{N}}}

\newcommand{\GL}{{\rm GL}}
\newcommand{\JSJ}{{\rm JSJ}}

\font\titsc=cmcsc10 scaled 1200

\newcommand{\ptwoirred}{$\matP^2$-irreducible}

\newcommand{\matr} [4] {\tiny{\left(\begin{array}{@{}c@{\ }c@{}} #1 & #2 \\ #3 & #4 \\ \end{array} \right)}}
\newcommand{\tr} {{\rm tr\,}}

\newcommand{\timtil}{\begin{picture}(12,12)
\put(2,0){$\times$}\put(2,4.5){$\sim$}\end{picture}}
\newcommand{\timtiltil}{\begin{picture}(12,12)
\put(2,0){$\times$}\put(2,4){$\sim$}\put(2,5.5){$\sim$}\end{picture}}

\newcommand{\finedimo}{{\hfill\hbox{$\square$}\vspace{2pt}}}
\newcommand{\dimo}[1]{\vspace{1pt}\noindent{\it Proof of} {\hspace{2pt}}\ref{#1}.}

\author{Gennaro \titsc{Amendola} 
\and Bruno \titsc{Martelli}}

\title{Non-orientable 3-manifolds of complexity up to $7$}

\begin{document}

\maketitle

\begin{abstract}
  We classify all closed non-orientable \ptwoirred\ 3-manifolds with complexity up to $7$,
  fixing two mistakes in our previous complexity-up-to-$6$ classification. 
  We show that there is no such manifold with complexity less than $6$, 
  five with complexity $6$
  (the four flat ones and the filling of the Gieseking manifold, which is of type \Sol),
  and three with complexity $7$ (one manifold of type \Sol, and the two
  manifolds of type $\matH^2\times \matR$ with smallest base orbifolds).\\[4pt]
\end{abstract}

\section*{Introduction}
It has been experimented in various contexts that 
non-orientable 3-manifolds are much more sporadic
than orientable ones. First of all, among the 8 three-dimensional geometries, only
$5$ have non-orientable representatives. Then, among cusped hyperbolic manifolds
of complexity up to $7$, only $1260$ of $6075$ are non-orientable, as shown in the
Callahan-Hildebrand-Weeks census~\cite{CaHiWe}. Here we show that, among
\emph{closed} \ptwoirred\ manifolds of complexity up to $7$, only $8$ of $318$
are non-orientable.

The complexity we refer to is the one defined by Matveev~\cite{Mat88, Mat90}.
As shown in~\cite{MaPe:nonori}, the complexity $c(M)$ of a closed \ptwoirred\ $M$ distinct
from $S^3, \matRP^3, L_{3,1}$ equals the
minimum number of tetrahedra needed to triangulate $M$.
Closed non-orientable \ptwoirred\ manifolds of complexity up to $6$ 
were classified in~\cite{AmMa} using only
theoretical arguments. The arguments were correct, except for
two mistakes in recognizing the geometries of the
resulting manifolds:
we fix them at the end of Section~\ref{main:section}.
(Namely, it is not true that all manifolds 
with complexity $c=6$ are flat,
and that there is one non-geometric manifold with $c=7$, as asserted in~\cite{AmMa}.)

The main result of this paper, stated in Theorem~\ref{main:teo} below, 
is the classification of all closed non-orientable \ptwoirred\
manifolds with complexity $c\leqslant 7$.
The contribution of this result to the census of all manifolds with
$c\leqslant 9$ is summarized in Table~\ref{tabella}.
Theorem~\ref{main:teo} is stated and proved in Section~\ref{main:section}. 
The proof of a lemma is deferred to
Section~\ref{lemma:section} and some facts on $I$-bundles over
surfaces and on $\matH^2\times\matR$-manifolds are collected in the
Appendix.

\begin{table}[tp]
  \begin{center}
    \begin{tabular}{r|c|c|c|c|c|c|c|c|c|c}
      & $0$ & $1$ & $2$ & $3$ & $4$ & $5$ & $6$ & $7$ & $8$ & $9$ \\
      
      \hline\hline

      \multicolumn{11}{c}{orientable\phantom{\Big|}} \\

      \hline
      \hline
      
      \phantom{\Big|} 
      lens &
      $3$ & $2$ & $3$ & $6$ & $10$ & $20$ & $36$ & $72$ & $136$ & $272$ \\

      \hline

      other elliptic &
      \multicolumn{1}{c}{\phantom{\Big|}} & & $1$ & $1$ & $4$ & $11$ & $25$ & $45$ & $78$ & $142$ \\

      \hline
      
      flat &
      \multicolumn{5}{c}{\phantom{\Big|}} & & $6$ & \multicolumn{3}{c}{} \\
      
      \hline
      
      Nil &
      \multicolumn{5}{c}{\phantom{\Big|}} & & $7$ & $10$ & $14$ & $15$ \\

      \hline

      $\widetilde{{\rm SL}_2\matR}$ &
      \multicolumn{6}{c}{\phantom{$\Big|^|$}} & & $39$ & $162$ & $514$ \\
      
      \hline

      \Sol &
      \multicolumn{6}{c}{\phantom{\Big|}} & & $5$ & $9$ & $23$ \\

      \hline

      $\matH^2\times\matR$ &
      \multicolumn{7}{c}{\phantom{\Big|}} & & $2$ & \\
      
      \hline 

      hyperbolic &
      \multicolumn{8}{c}{\phantom{\Big|}} & & $4$ \\
      
      \hline

      non-trivial \JSJ &
      \multicolumn{6}{c}{\phantom{\Big|}} & & $4$ & $35$ & $185$ \\
      
      \hline
      \hline 

      total orientable \phantom{\Big|} &
      $3$ & $2$ & $4$ & $7$ & $14$ & $31$ & $74$ & $175$ & $436$ & $1155$ \\

      \hline
      \hline

      \multicolumn{11}{c}{{\bf non-orientable}\phantom{\Big|}} \\

      \hline
      \hline

      flat &
      \multicolumn{5}{c}{\phantom{\Big|}} & & $\bf{4}$ & \multicolumn{3}{c}{} \\

      \hline

      $\matH^2\times\matR$ &
      \multicolumn{6}{c}{\phantom{\Big|}} & & $\bf{2}$ & \multicolumn{2}{c}{?} \\ 

      \hline

      \Sol &
      \multicolumn{5}{c}{\phantom{\Big|}} & & $\bf{1}$ & $\bf{1}$ & \multicolumn{2}{c}{?} \\

      \hline
      \hline 

      total non-orientable &
      \multicolumn{5}{c}{\phantom{\Big|}} & & $\bf{5}$ & $\bf{3}$ & \multicolumn{2}{c}{?} \\

      \hline

    \end{tabular}
  \end{center} 
  \caption{The number of \ptwoirred\ manifolds of
    given complexity (up to $9$) and geometry (empty boxes contain
    $0$).}
  \label{tabella}
\end{table}

Theorem~\ref{main:teo} has been proved independently by
Burton~\cite{Bu} using the computer program \emph{Regina}~\cite{Re}. 
More than that, Burton has
classified all minimal triangulations with at most $7$ tetrahedra.

\paragraph{Acknowledgements} 
We thank Maria Rita Casali and Sergej Matveev for describing to us the
\Sol-manifold missing in our previous list~\cite{AmMa}.

\section{Main statement}\label{main:section}

We recall that there are 8 important 3-dimensional geometries, six of them 
concerning Seifert manifolds. 
The geometry of a Seifert manifold is determined by two invariants of any of its fibrations,
namely the Euler characteristic $\chiorb$ of the base orbifold and the Euler number
$e$ of the fibration, according to Table~\ref{tabellina}. The two non-Seifert
geometries are the hyperbolic and the \Sol\ ones.
We refer to~\cite{Sco} for definitions.

\begin{table}[tp]
  \begin{center}
    \begin{tabular}{c|ccc} 
      \phantom{\Big|} & $\chiorb>0$ & $\chiorb=0$ & $\chiorb<0$ \\ \hline
      \phantom{\Big|} $e=0$ & $\matS^2\times\matR$ & $\matE^3$ & $\matH^2\times\matR$ \\
      \phantom{\Big|} $e\neq 0$ & $\matS^3$ & Nil & $\widetilde{{\rm SL}_2\matR}$ \\ 
    \end{tabular}
  \end{center}
  \caption{The six Seifert geometries.}
  \label{tabellina}
\end{table} 

The complete list of closed \emph{orientable} irreducible manifolds of complexity $c\leqslant 9$ is
available from~\cite{weblist} and summarized in the first half of Table~\ref{tabella}.
The second half of Table~\ref{tabella} is recovered from the next result.

\begin{teo} \label{main:teo}
There are no non-orientable \ptwoirred\ manifolds with complexity $c\leqslant 5$.
There are $5$ such manifolds with $c=6$: they are the $4$ flat ones and
the torus bundle (of type \Sol) with monodromy $\matr 1110$.
There are $3$ such manifolds with $c=7$:
they are the torus bundle (of type \Sol) with monodromy $\matr 2110$, and the
two non-orientable Seifert manifolds (of type $\matH^2\times\matR$) defined by
$$\left(\matRP^2,(2,1),(3,1)\right) \quad {\rm and} \quad \left(\bar D,(2,1),(3,1)\right).$$
\end{teo}

Concerning the statement of Theorem~\ref{main:teo}, we denote by $\bar D$ the orbifold given by the 
disc with mirrored boundary.
Moreover, well-definition of the two non-orientable Seifert manifolds of
type $\matH^2\times\matR$ is proved in~\cite[pages 15 and 90]{Or}. Using the notations of~\cite{Or}, 
the two manifolds are respectively
$$\left\{0;({\rm o}, 0,0,1);(2,1),(3,1)\right\} \qquad {\rm and} 
\qquad \left\{0;({\rm n}_1,1);(2,1),(3,1)\right\}.$$

\begin{rem} {\em
The non-orientable \ptwoirred\ manifolds with $c\leqslant 7$ are the ``simplest'' ones in
each geometry:
the Gieseking manifold (the cusped hyperbolic manifold with smallest volume $1.0149\ldots$~\cite{Ad} and
smallest complexity $1$~\cite{CaHiWe}) is the punctured torus bundle with monodromy $\matr 1110$. Therefore the 
\Sol-manifold with $c=6$ is the (unique) filling (with a solid Klein bottle) of the Gieseking
manifold. The two \Sol-manifolds with $c\leqslant 7$ are the only torus bundles over $S^1$ whose
monodromy $A$ is hyperbolic with $|\tr A|\leqslant 2$, see Proposition~\ref{Sol:cor} in the Appendix.
The two $\matH^2\times\matR$-manifolds with $c=7$ have the smallest possible base hyperbolic orbifold,
having volume $-2\pi\chi^{\rm orb} = \pi/3$, see
Proposition~\ref{small:Seifert:prop} in the Appendix.
}
\end{rem}

The main ingredient in the proof of Theorem~\ref{main:teo} is the known 
list of all \emph{orientable} irreducible manifolds
with $c\leqslant 9$, available from~\cite{weblist}, and the
following lemma, which holds for manifolds of any complexity. 

\begin{lemma} \label{main:lemma}
Let $M$ be a closed non-orientable \ptwoirred\ manifold, and $\widetilde M$ be its
orientable double covering. We have $c(\widetilde M)\leqslant 2\cdot c(M)-5$.
\end{lemma}

The proof of Lemma~\ref{main:lemma} is deferred to Section~\ref{lemma:section}. The rest
of this section is devoted to the proof of Theorem~\ref{main:teo}.

\paragraph{Geometric decomposition}
We denote since now by $T$ and $K$ respectively the torus and the Klein bottle.
Let $M$ be a closed \ptwoirred\ manifold.
We recall that $M$ has a unique \emph{geometric decomposition} along
embedded tori and Klein bottles, defined as follows:
take the set of tori and Klein bottles of the \JSJ\ decomposition, and substitute each 
element of this set bounding an $I$-bundle over $T$ or $K$ with the core $T$ or $K$. 
In contrast to the
\JSJ, the geometric decomposition has two nice properties: it decomposes
$M$ into \emph{blocks} with finite volume,
and it remains geometric when lifted to finite coverings of $M$.
See Corollaries~\ref{negative:base:cor} and~\ref{decomposition:lifts:cor} in
the Appendix.

\paragraph{Seifert blocks} 
Let now $M$ be closed non-orientable and $\widetilde M$ be its orientable double-covering.
As we said, the geometric decomposition of $M$ lifts to the one of $\widetilde M$. 
Let $N$ be a block of the decomposition of $M$. Its pre-image in $\widetilde M$ is amphichiral,
\emph{i.e.}~it admits an orientation-reversing involution. Let us fix an orientation on $\widetilde M$.
The pre-image of $N$ consists of two blocks or one block, depending on
whether $N$ is orientable or not.
If $N$ is Seifert, its pre-image has Euler number zero~\cite{Neu}. (If it
consists of two blocks $\widetilde N_1$ and $\widetilde N_2$,
we mean that $e(\widetilde N_1)=-e(\widetilde N_2)$).
In particular, if the whole $M$ is itself Seifert, both
$M$ and $\widetilde M$ are
either flat or of type $\matH^2\times \matR$.

\paragraph{Orientable coverings of small complexity}
The following result, together with Lemma~\ref{main:lemma} and Proposition~\ref{upper:bounds:prop}
below, will easily imply Theorem~\ref{main:teo}.
\begin{prop} \label{main:prop}
Let $M$ be a closed non-orientable \ptwoirred\ manifold. If 
its orientable double-covering $\widetilde M$ has complexity $c(\widetilde M)\leqslant 9$,
then one of the following occurs:
\begin{itemize}
\item $c(\widetilde M)=6$, and both $\widetilde M$ and $M$ are flat;
\item $c(\widetilde M)=7$, and both $\widetilde M$ and $M$ are \Sol\
  torus bundles over $S^1$, with monodromies $\matr 2111$ and
  $\matr1110$;
\item $c(\widetilde M)=8$, both $\widetilde M$ and $M$ are of type $\matH^2\times\matR$, and 
\begin{itemize}
\item
$\widetilde M$ is $\big(S^2,(2,1),(3,1),(2,-1),(3,-1)\big)$, 
\item
$M$ is either $\big(\matRP^2,(2,1),(3,1)\big)$ or $\big(\bar D,(2,1),(3,1)\big);$
\end{itemize}
\item $c(\widetilde M)=9$, and both $\widetilde M$ and $M$ are \Sol\
  torus bundles over $S^1$, 
with monodromies $\matr 5221$ and $\matr2110$.
\end{itemize}
\end{prop}
\begin{proof}
We denote by $D$, $A$, and $S$ respectively the disc, the annulus, and
the M\"obius strip.
Since $M$ is \ptwoirred, the orientable double-covering $\widetilde M$
is irreducible.
Now, an orientable irreducible manifold with complexity $c\leqslant 9$
has one of the following geometric decompositions~\cite{weblist}:
\newcounter{listi}
\begin{list}{{\rm (\roman{listi})}}{\usecounter{listi}\setlength{\labelwidth}{1cm}}
\item \label{SS:list}
it is itself Seifert or \Sol;
\item \label{DD:list}
it decomposes along one $T$ into two Seifert blocks, each fibering over
$D$ with two singular fibers; 
\item \label{D:list}
it decomposes along one $K$ into one Seifert block fibering
over $D$ with two singular fibers;
\item \label{c9:list}
it decomposes along one or two $K$'s into one Seifert block, which fibers either over $D$
with 3 singular fibers of type (2,1), or over $S$ or $A$ with one singular fiber of type (2,1);
\item \label{H:list}
it is one of the $4$ smallest hyperbolic manifolds known.
\end{list}

Cases~(ii-v) occur only for $c\geqslant 7$.
Note that the \JSJ\ decomposition used in~\cite{weblist} should be translated into the geometric
one by replacing each block of type $(D,(2,1),(2,1))$ with a $K$, thus getting 
cases~(iii) and~(iv). Cases~(iv-v) only occur for $c=9$.

Suppose $\widetilde M$ is of type~(i). If it is Seifert, since $e(\widetilde M)=0$, 
it is either flat or of type $\matH^2\times \matR$.
If $\widetilde M$ is flat, we are done. Suppose it is of type $\matH^2\times \matR$.
There are only two such manifolds in the list: they both have $c=8$ and they are
$$\big(S^2,(2,1),(3,1),(2,-1),(3,-1)\big) \quad {\rm and } \quad 
\big(\matRP^2,(3,1),(3,-1)\big).$$
In both cases we have $\chiorb(\widetilde M) = -1/3$. Therefore
$\chiorb(M)=-1/6$.
Now, Proposition~\ref{small:Seifert:prop}, proved
in the Appendix, shows that there are two possible $M$'s.
They have the same double cover, as required.

If $\widetilde M$ is \Sol, then $M$ is \Sol\ too. By Corollary~\ref{Sol:cor}, proved in the
Appendix, both $\widetilde M$ and $M$ are torus bundles over $S^1$ with some monodromies $A^2$ and $A$.
From linear algebra we get $\tr(A^2) = (\tr A)^2 - 2\det A = (\tr A)^2 + 2$. 
The orientable manifolds
in the list~\cite{weblist} satisfy $3\leqslant |\tr(A^2)| \leqslant
8$.
The only possible values for
$\tr (A^2)$ are then $3$ and $6$ (namely, $A^2$ is $\matr 2110$ or
$\matr 5221$),
so $|\tr A|\in\{1,2\}$. By Corollary~\ref{Sol:cor},
there is only one non-orientable manifold for each such value of $|\tr A|$, hence
we get two manifolds, with monodromies $\matr1110$ and $\matr2110$.

We are left to prove that $\widetilde M$ cannot be of types~(ii-v). Suppose it is of type~(ii).
Then $\widetilde M$ is the union of two Seifert manifolds
$\widetilde N_i = \big(D,(p_i,q_i),(r_i,s_i)\big)$, with $r_i>2$ and $i\in\{1,2\}$.
By what said above, the geometric decomposition of $M$ 
consists of either one block homeomorphic to both $\widetilde N_1$ and $\widetilde N_2$, or of two blocks $N_1$ and $N_2$,
with $\widetilde N_i$ covering $N_i$. If the first possibility holds, there is an involution 
$\tau:\widetilde M\to\widetilde M$ exchanging $\widetilde N_1$ and $\widetilde N_2$ and giving $M$ as a quotient.
That $\tau$ restricts to an orientation-preserving order-2
involution on the torus $T$ separating $\widetilde N_1$ and
$\widetilde N_2$.
Therefore $\tau$ acts like $\pm I$ on
$H_1(T)$, thus preserving simple closed curves (up to isotopy). On the other side, 
$\tau$ sends a fiber of $\widetilde N_1$ to a fiber of $\widetilde N_2$, but these fibers
give non-isotopic curves on $T$, and we get a contradiction.
If the second possibility holds, we have $e(\widetilde N_1)=e(\widetilde N_2)=0$,
hence $p_i=r_i>2$ for $i=1,2$.
But no manifold with $c\leqslant 9$ in the list~\cite{weblist} has
these parameters.

If $\widetilde M$ is of type~(iii), it is decomposed along a single $K$.
But a manifold whose decomposition contains an odd number of $K$'s is not the double covering
of a non-orientable one, see Corollary~\ref{even:cor}.
Finally, $\widetilde M$ cannot be of type~(iv) because the unique Seifert block has Euler number
$e=1/2\neq 0$. And it cannot be of type~(v), because the deck involution would be an isometry, 
but there is no orientation-reversing isometry of the $4$ smallest
closed hyperbolic manifolds known giving a manifold~\cite{HoWe}.
\end{proof}

\paragraph{Spines and complexity}
We briefly recall some definitions from~\cite{Mat90}. 
A compact 2-dimensional
polyhedron $P$ is \emph{simple} if the link of every point in $P$ is contained in
the 1-skeleton $K$ of the tetrahedron.
A point having the whole of $K$ as a link is 
called a \emph{vertex}. The set $V(P)$ of the vertices of $P$ consists of isolated points, so
it is finite. A compact polyhedron $P\subset M$ is a \emph{spine} of the closed manifold $M$
if $M\setminus P$ is an open ball. The \emph{complexity} $c(M)$ of
a closed 3-manifold $M$ is then defined as the minimal number of vertices of
a simple spine of $M$. It turns out~\cite{Mat90, MaPe:nonori} that
if $M$ is \ptwoirred\ and distinct from $S^3, \matRP^3, L_{3,1}$ then it has a \emph{minimal} spine 
(\emph{i.e.}~a spine with $c(M)$ vertices)
which is \emph{special}. A spine $P$ is special when it is the 2-skeleton of the dual of a 1-vertex triangulation of $M$.
Its singular set $S(P)$ is a connected 4-valent graph.

\paragraph{Manifolds with marked boundary}
We now recall some definitions from~\cite{MaPe:nonori}, which we will use to prove
Proposition~\ref{upper:bounds:prop} below.
A \emph{$\theta$-graph} in the torus $T$ is a trivalent graph $\theta\subset T$ whose complement
in $T$ is an open disc.
Let $M$ be a connected (possibly non-orientable) compact 3-manifold with (possibly empty)
boundary consisting of tori. By associating to each component of
$\partial M$ a $\theta$-graph, we get a
{\em manifold with marked boundary} (one can also define markings on Klein bottles, 
see~\cite{MaPe:nonori}). A simple polyhedron $P\subset M$ which intersects
$\partial M$ in the union of the markings and such that $M\setminus (P\cup\partial M)$
is an open ball is a \emph{skeleton} for the marked $M$. When $M$ is closed, a skeleton
is just a spine.

Given two marked $M, M'$ and a homeomorphism $\psi$ between one component of $\partial M$ and
one of $\partial M'$ which preserves the markings, one can glue $M$ and $M'$ and get a new
marked manifold $N$, which is called an \emph{assembling} of $M$ and $M'$. 
Two skeleta $P$, $P'$ of $M$, $M'$ glue via $\psi$ to a skeleton $Q$
of $N$. Spines of plenty of manifolds
can be constructed in this way, and by controlling the number of their vertices one gets
many strict upper bounds for complexity~\cite{MaPe:geometric}.
Here, we need the following one.

\begin{prop} \label{upper:bounds:prop}
Every flat non-orientable manifold has complexity $c\leqslant 6$. The torus bundles
with monodromy $\matr 1110$ and $\matr 2110$ have respectively $c\leqslant 6$ and $c\leqslant 7$.
The closed non-orientable manifolds
$$\big(\matRP^2,(2,1),(3,1)\big) \quad {\rm and} \quad \big(\bar D,(2,1),(3,1)\big)$$
have complexity $c\leqslant 7$.
\end{prop}
\begin{proof}
Spines of flat manifolds with $6$ vertices are constructed in~\cite[Section 3]{AmMa}.
The upper bound for torus bundles $M$ with monodromy $A\in\GL_2(\matZ)$ 
given in~\cite{MaPe:geometric} works also in the non-orientable
case, and it gives $c(M)\leqslant\max\{||A||+5, 6\}$. From the definition of the norm $||A||$
in~\cite{MaPe:geometric} one easily gets $||\matr k110||\leqslant k$ for $k>0$, as required. 

Finally, by Proposition~\ref{I-bundles:prop} proved in the Appendix,
each of the two Seifert manifolds is the result of gluing $N=\big(D,(2,1),(3,1)\big)$
to $T\timtil I$ with an appropriate map. A skeleton for $N$ such that
the marking $\theta$ contains a loop $\gamma$ isotopic to the fiber is
constructed in~\cite[Page 170]{AmMa} and also shown in Fig.~\ref{TtimtilI:fig}-left: it has $4$ vertices.
A skeleton with 3 vertices of a marked $T\timtil I$ is
shown in Fig.~\ref{TtimtilI:fig}-right.
The marking $\theta'$ contains two loops $\gamma'_1, \gamma'_2$
isotopic to the two distinct fibrations of $T\timtil I$.
Therefore it is possible to assemble $N$ and $T\timtil I$
sending $\gamma$ either to $\gamma'_1$ or to $\gamma'_2$, and the two assemblings give the two
Seifert manifolds above.
\begin{figure}[tp]
  \begin{center}
    \mettifig{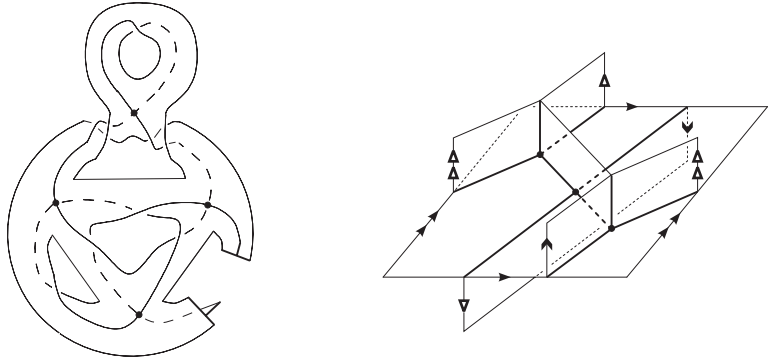}
    \caption{A skeleton for $\big(D,(2,1),(3,1)\big)$ with 4 vertices (left) and one for $T\widetilde\times I$ with 3 vertices (right).}
    \label{TtimtilI:fig}
  \end{center}
\end{figure}
\end{proof}

We can finally prove Theorem~\ref{main:teo}.

\dimo{main:teo}
Let $M$ be a closed non-orientable \ptwoirred\ manifold. Then its orientable double cover
$\widetilde M$ is irreducible. Lemma~\ref{main:lemma} gives
$c(\widetilde M)\leqslant 2c(M)-5$. If $c(M)\leqslant 5$, we get $c(\widetilde M)\leqslant 10-5=5$,
which is impossible by Proposition~\ref{main:prop}.
If $c(M)=6$, we get $c(\widetilde M)\leqslant 12-5=7$. By Proposition~\ref{main:prop},
the manifold $M$ is either flat or a torus bundle with monodromy
  $\matr 1110$.
Now, note that Lemma~\ref{main:lemma} does not guarantee the converse,
namely that all such manifolds have $c=6$.
But Proposition~\ref{upper:bounds:prop} gives
$c\leqslant 6$ on them, hence $c=6$, as required.

If $c(M)=7$, we get $c(\widetilde M)\leqslant 14-5=9$.
By Proposition~\ref{main:prop} (and by what said above),
the manifold $M$ is either a torus bundle with monodromy $\matr 2110$ or 
one of the two Seifert manifolds of type $\matH^2\times\matR$.
Again, each such manifold has $c\leqslant 7$ by Proposition~\ref{upper:bounds:prop}, and
we are done.
\finedimo

\paragraph{Errata}
We now fix two mistakes present in our previous paper~\cite{AmMa}.
First, it is stated there that closed non-orientable \ptwoirred\
manifolds of complexity $6$ are flat, 
because the torus bundle with monodromy $\matr 1110$ was not recognized
as a \Sol\ manifold in~\cite[Section 3, page 169]{AmMa}.
Second, an example of closed non-orientable \ptwoirred\ manifold $M$
with complexity $7$ having non-trivial \JSJ\ is shown in~\cite[Section
  3, page 170]{AmMa}.
It would consist of two Seifert
blocks, one being the non-orientable $I$-bundle $T\timtil I$ over $T$. As shown
in Proposition~\ref{I-bundles:prop} here, every gluing of $T\timtil I$ with a Seifert manifold
is itself Seifert. Therefore the \JSJ\ of $M$ is indeed trivial (actually, $M$ is one of
the two manifolds of type $\matH^2\times\matR$ having $c=7$).

\section{Stiefel-Whitney surfaces} \label{lemma:section}

This section is devoted to the proof of Lemma~\ref{main:lemma}.
We start with some preliminary results.

\paragraph{Stiefel-Whitney surfaces}
A closed non-orientable manifold $M$ has a non-trivial first
Stiefel-Whitney class $w_1\in H^1(M;\matZ/_{2\matZ})$.
A surface $\Sigma\subset M$ which is Poincar\'e dual to $w_1$ is
usually called a {\em Stiefel-Whitney
surface}~\cite{HeLaNu}.
It has odd intersection with a transverse loop $\gamma$ if and only if $\gamma$
is orientation-reversing. It is easy to prove that $\Sigma$ is orientable.
Note that there are infinitely many non-isotopic Stiefel-Whitney
surfaces in $M$.

We will now show that, fixed a special spine $P$ of $M$, there is
exactly one Stiefel-Whitney surface contained in $P$.
The embedding $P\subset M$ induces an isomorphism
$H_2(P;\matZ/_{2\matZ}) \cong H_2(M;\matZ/_{2\matZ})$.
Using cellular homology, a representative for a cycle in
$H_2(P;\matZ/_{2\matZ})$ is a subpolyhedron consisting of some faces,
an even number of them (whence 0 or 2) incident to each edge of $P$.
Such a subpolyhedron is a surface near the edges it contains, and it
is also a surface near the vertices (because the link of a vertex
does not contain two disjoint circles).
Thus every $\matZ/_{2\matZ}$-homology class is represented by a unique
surface in $P$: in particular there is a unique Stiefel-Whitney
surface $\Sigma(P)$ inside $P$.

\begin{rem} \label{lift:rem} {\em
Let $P$ be a special spine of a non-orientable closed $M$, and $\Sigma = \Sigma(P)$ be
the Stiefel-Whitney surface contained in $P$.
Let $\widetilde \Sigma\subset \widetilde P\subset \widetilde M$ be the pre-images of 
$\Sigma\subset P\subset M$ in the orientable double-cover $\widetilde M$. Both
$\Sigma$ and $\widetilde\Sigma$ are orientable.
Here $\widetilde M\setminus\widetilde P$ consists of two balls, and $\widetilde\Sigma$
consists precisely of all faces of $\widetilde P$ that are adjacent to both these balls. 
Making a hole on one face contained in $\widetilde\Sigma$ 
one gets a $\widetilde P'$ whose
complement in $\widetilde M$ is a single ball, \emph{i.e.}~a simple spine of $\widetilde M$.
}
\end{rem}

By Remark~\ref{lift:rem}, if $P$ is a minimal spine of $M$ with $v$ vertices, 
there is a spine $\widetilde P'$ for
$\widetilde M$ with $2v$ vertices. This gives $c(\widetilde M)\leqslant M$. 
But $\widetilde P'$ has a hole in a face 
$F\subset\widetilde\Sigma$,
which can be enlarged with a collapse, eventually deleting the whole $F$ and killing 
all the vertices adjacent to $F$.
The number of such vertices
killed depends on the choice of $F$ in $\widetilde\Sigma$. The rest
of this section is devoted to the proof that there is one face 
$F\subset \widetilde \Sigma$ incident to at least
$5$ distinct vertices. Using such an $F$, we get 
a simple spine for $\widetilde M$ with $2v-5$ vertices at most, hence proving Lemma~\ref{main:lemma}.

\paragraph{Length of a face}
Let $F$ be a face of a special spine $P$.
We denote by $\lgh{F}$ the number of vertices of $P$ adjacent to $F$,
counted with multiplicity.

\begin{lemma}\label{length:lemma}
  Let $P$ be a special spine of a closed non-orientable \ptwoirred\ $M$.
  Let $\widetilde\Sigma$ and $\widetilde P$ be the pre-images of
  $\Sigma=\Sigma(P)$ and $P$ in $\widetilde M$.
  There exists a face $F\subset\widetilde\Sigma$ with
  $\lgh{F}\geqslant 5$.
\end{lemma}
\begin{proof}
  The average value of $\lgh{F}$ on the faces in
  $\widetilde \Sigma$ is $s/f$, where $f$ is the number of faces of
  $\widetilde P$ contained in $\widetilde\Sigma$ and
  $s=\sum_{ F\subset\widetilde\Sigma}\lgh{ F}$.
  We prove that $s/f>4$, thus getting a face $F$ with
  $\lgh{ F}\geqslant 5$.
  Let $n_3$ be the number of pairs of 3-valent vertices of
  $\widetilde G=S(\widetilde P)\cap\widetilde\Sigma$ and $n_4$ be the
  number of 4-valent ones.
  The graph $\widetilde G$ has $2n_3+n_4$ vertices and
  $\frac{3(2n_3)+4n_4}{2} = 3n_3+2n_4$ edges, so
  $\chi(\widetilde\Sigma) = (2n_3+n_4) - (3n_3+2n_4) + f$.
  Hence $f$ is equal to $\chi(\widetilde\Sigma) + n_3 + n_4$.
  Moreover, the sum $s$ is equal to $6n_3+4n_4$, so the average value
  $s/f$ is 
  $$\frac sf = \frac{6n_3+4n_4}{\chi(\Sigma)+n_3+n_4}.$$
  Now, $\Sigma$ is orientable and non-separating (because $M\setminus
  P$ is a ball), and $M$ is \ptwoirred, so we get
  $\chi(\widetilde\Sigma)=2\chi(\Sigma)\leqslant 0$.
  Therefore, we have $s/f\geqslant 4$, with $s/f=4$ if and only if
  $\chi(\widetilde\Sigma)=0$ and $n_3=0$.
  But in the last case $\widetilde\Sigma$ would be a torus, and
  $\widetilde P$ would be the union of a torus with
  two discs, hence $\widetilde M$ would have genus $\leqslant 1$.
  This gives a contradiction, since both $\widetilde M$ and $M$ would
  be elliptic or of type $\matS^2\times\matR$.
\end{proof}

\paragraph{The polyhedron $\widetilde P$ near $\widetilde\Sigma$}
Let $P$ be a minimal spine of a closed non-orientable \ptwoirred\
$M$.
Lemma~\ref{length:lemma} guarantees the existence of a face of length at least $5$ in the pre-image 
$\widetilde\Sigma$ of $\Sigma(P)$.
Unfortunately, such a face might not be embedded, and hence might be incident to less than $5$ vertices.

\begin{figure}[tp]
\begin{center}
  \mettifig{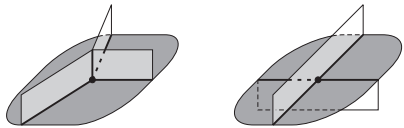}
  \caption{The spine $\widetilde P$ near a 3-valent (left) and a
    4-valent (right) vertex of $S(\widetilde P)\cap\widetilde\Sigma$.}
  \label{near_vert:fig}
\end{center}
\end{figure}

We now study the properties of non-embedded faces in $\widetilde\Sigma$. To do that,
we need to draw the spine $\widetilde P$ near the surface
$\widetilde\Sigma$.
The graph $\widetilde G=S(\widetilde P)\cap\widetilde\Sigma$ has
vertices with valence 3 and 4, and $\widetilde P$ appears near them as
shown in Fig.~\ref{near_vert:fig}.
By Remark~\ref{lift:rem}, the surface $\widetilde\Sigma$ is
orientable, then we can choose a transverse orientation and give each
edge $e$ of $\widetilde G$ a black or grey color, depending
on whether $\widetilde P$ locally lies on the positive or on the
negative side of $\widetilde\Sigma$ near $e$.
A $3$-valent vertex is adjacent to edges with the same color, and a
$4$-valent vertex is adjacent to two opposite grey edges and two
opposite black ones.
Now, the regular neighborhood $\calN(\widetilde G)$ of $\widetilde G$ in
$\widetilde P$ can be immersed into $\matR^3$ so that
$\widetilde\Sigma\cap \calN(\widetilde G)$ is ``horizontal''.
The polyhedron $\calN(\widetilde G)$ is determined unambiguously by that
immersed graph and also the regular neighborhood $\calN(\widetilde\Sigma)$
of $\widetilde\Sigma$ in $\widetilde P$ is, because it is obtained
from $\calN(\widetilde G)$ by adding discs to the ``horizontal'' $S^1$'s
in the boundary of $\calN(\widetilde G)$.

\begin{lemma}\label{no_loop:lemma}
  Let $P$ be a minimal spine of a closed non-orientable \ptwoirred\
  $M$.
  Then each edge in $G=S(P)\cap\Sigma$ and $\widetilde
  G=S(\widetilde P)\cap\widetilde\Sigma$ has different endpoints.
\end{lemma}
\begin{proof}
  Suppose that there exists an edge $e$ of $G$
  joining a vertex $v$ of $P$ to itself. The closure $\bar e$ of $e$ is then a loop in $\Sigma$.
  We prove that $\bar e$ bounds an embedded disc $D\subset P$ with $\lgh D=1$,
  which is absurd because $P$ is minimal~\cite{Mat90}.
  If $v$ is 4-valent and the two germs of $e$ near $v$ are opposite, the
  disc $D$ lies in $P\setminus\Sigma$.
  If not, the two germs are ``consecutive'' near $v$. Since $\Sigma$ is orientable, the neighborhood of $\bar e$
  in $\Sigma$ is an annulus, hence $\bar e$
  bounds a disc $D\subset\Sigma$. The same result for $\widetilde G$ follows.
\end{proof}

\begin{lemma}\label{no_small_faces:lemma}
  Let $P$ be a minimal spine of a closed non-orientable \ptwoirred\
  $M$.
  There is no embedded face $F$ of $P$ with $\partial
  F\subset\widetilde\Sigma$ and $\lgh F\leqslant 3$.
  If moreover $\Sigma(P)$ contains the minimum number
  of vertices (among minimal spines $P$ of $M$), then there is no
  embedded square $F\subset\widetilde\Sigma\subset\widetilde P$ with the following
  shape:\\[10pt]
  \centerline{\mettifig{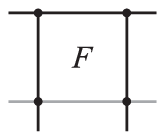}}
\end{lemma}
\begin{proof}
  By Lemma~\ref{no_loop:lemma}, two consecutive vertices in $\partial F\subset\widetilde\Sigma$ project to
  two distinct vertices in $\Sigma$. Therefore, if $\lgh{F}\leqslant 3$ all vertices of $\partial F$
  project to distinct vertices of $P$, hence $F$ projects to an embedded face
  with ${\rm lgh}\leqslant 3$, in contrast to
  minimality of $P$~\cite{Mat90}.
  Suppose now $F$ is a square as above.
  Opposite vertices of $\partial F$ have distinct valency, 
  hence they project to distinct vertices of $\Sigma$.
  Therefore, the projection of $F$ is an embedded square in $\Sigma$,
  and the move shown in Fig.~\ref{lgh_5_move:fig} transforms
  $P$ into another minimal $P'$, but with $\Sigma(P')$ containing one vertex less 
  than $\Sigma(P)$, a contradiction.
  \begin{figure}[tp]
  \begin{center}
    \mettifig{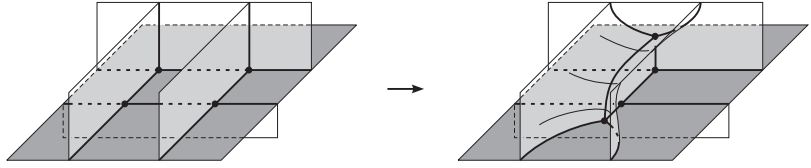}
    \caption{A move contradicting the assumption that the number of
      vertices of $P$ contained in $\Sigma$ is minimal.}
    \label{lgh_5_move:fig}  
  \end{center}
  \end{figure}
\end{proof}

\begin{lemma}\label{no_self_adj_edge:lemma}
  Let $P$ be a minimal spine of a closed non-orientable \ptwoirred\
  $M$.
  Then each face
  $ F\subset\widetilde\Sigma\subset\widetilde P$ is not
  incident twice to an edge of
  $\widetilde G=S(\widetilde P)\cap\widetilde\Sigma$.
\end{lemma}
\begin{proof}
  Suppose by contradiction that there is a face
  $F\subset\widetilde\Sigma\subset\widetilde P$ incident
  twice to an edge $e$ of $\widetilde G$.
  Inside the closure of $F$, there is a loop
  $\lambda$ intersecting $S(\widetilde P)$ transversely in
  one point of $ e$.
  The regular neighborhood of $\lambda$ in
  $\widetilde\Sigma$ is transversely orientable (because
  $\widetilde\Sigma$ is orientable), hence $\lambda$ bounds
  a disc $ D$ in $\widetilde M\setminus\widetilde P$ and
  intersects $S(\widetilde P)$ transversely in 1 point.
  If we project $\lambda$ to $\Sigma$, we get a loop
  $\lambda$ which bounds a disc (the projection of $D$) in
  $M\setminus P$ and intersects $S(P)$ transversely in
  1 point: this contradicts the minimality of $P$~\cite{Mat90}.
\end{proof}

We can finally prove Lemma~\ref{main:lemma}.

\dimo{main:lemma}
  Let $P$ be a minimal spine of $M$, such that $\Sigma=\Sigma(P)$
  contains the minimum possible
  number of vertices (among minimal spines of $M$). Let 
  $\widetilde\Sigma\subset\widetilde P\subset\widetilde M$ be the pre-images
  of $\Sigma\subset P\subset M$.
  By what said in Remark~\ref{lift:rem} and below, if we prove that
  $\widetilde\Sigma$ contains a face incident to $5$ \emph{distinct}
  vertices at least we are done.
  Let $F$ be a face of $P$ contained in $\widetilde\Sigma$
  such that $\lgh{F}$ is maximal.
  By Lemma~\ref{length:lemma}, we have
  $\lgh{F}\geqslant 5$.
  If $ F$ is embedded, we are done.
  If instead $F$ is not embedded, we will show that there
  are only a finite (small) number of configurations of
  $\widetilde\Sigma$ near $ F$, and for each case we
  will find a face incident to $5$ distinct vertices in $\widetilde\Sigma$
  (or get a contradiction).

  So, from now on, we can suppose $F$ is not embedded.
  As we said above, if $F$ is incident to at least $5$
  \emph{distinct} vertices, we are done.
  So we are left to deal with the case where $ F$ is incident to at
  most $4$ distinct vertices.
  By Lemma~\ref{no_self_adj_edge:lemma}, $ F$ cannot be
  incident twice to an edge of the graph $\widetilde G = 
  S(\widetilde P)\cap\widetilde\Sigma$, so
  $F$ can be incident only once to a 3-valent vertex and
  either once or twice to a 4-valent one of $\widetilde G$.
  Since $F$ is not embedded, it is incident twice to at
  least one 4-valent vertex.
  We conclude the proof with a case-by-case argument.

  If $\lgh{F}\geqslant 9$, since $F$ is incident
  to each vertex at most twice, our $F$ would be incident
  to 5 different vertices of $\widetilde P$, a contradiction.
  If instead $\lgh{F}=8$, our $F$ is incident to 4
  different 4-valent vertices twice, so $S(\widetilde P)=\widetilde
  G=\partial  F$ (because $S(\widetilde P)$ is connected)
  and $\widetilde P$ has $4$ vertices, hence $\widetilde M$ is
  elliptic~\cite{Mat90, weblist}, a contradiction.
  
\begin{figure}[tp]
    \begin{center}
    \mettifig{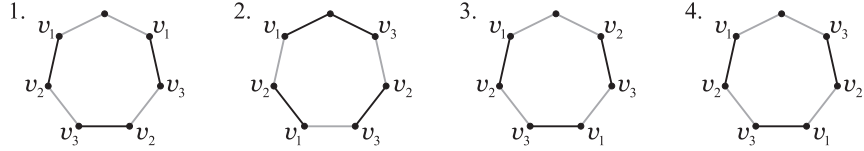}
    \caption{The four different configurations for the pairing of the
      vertices adjacent to $ F$, if $\lgh{ F}=7$.}
    \label{lgh_7_unfold:fig}
    \end{center}
  \end{figure}

  If $\lgh{F}=7$, our $F$ is incident to 4 different
  vertices: twice to 3 of them (which are 4-valent), and once to
  another one.
  If we consider the unfolded version of $F$, we have a
  heptagon with six vertices identified in pairs.
  Up to symmetry, there are 4 different configurations for the pairing
  of the vertices adjacent to $F$ (recall that
  Lemma~\ref{no_loop:lemma} forbids edges in $\partial  F$
  with coinciding endpoints).
  They are shown in Fig.~\ref{lgh_7_unfold:fig}.
  As we have done above, the black or grey color given to each edge
  depends on whether $\widetilde P$ locally lies on the positive or on
  the negative side of $\widetilde\Sigma$ near the edge.
  Recall that all the $ v_i$'s are 4-valent, so the two
  consecutive edges going out from a $ v_i$ have different
  colors.
  In each case, using orientability of $\widetilde\Sigma$, 
  one finds two edges in the boundary of the unfolded version of $F$
  that map to the same edge of $\widetilde P$, contradicting
  Lemma~\ref{no_self_adj_edge:lemma}.

\begin{figure}[tp]
    \begin{center}
    \mettifig{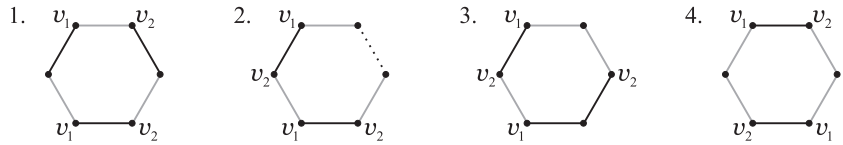}
    \caption{Four different configurations when $\lgh{ F}=6$.}
    \label{lgh_6_unfold:fig}
    \end{center}
  \end{figure}

\begin{figure}[tp]
    \begin{center}
    \mettifig{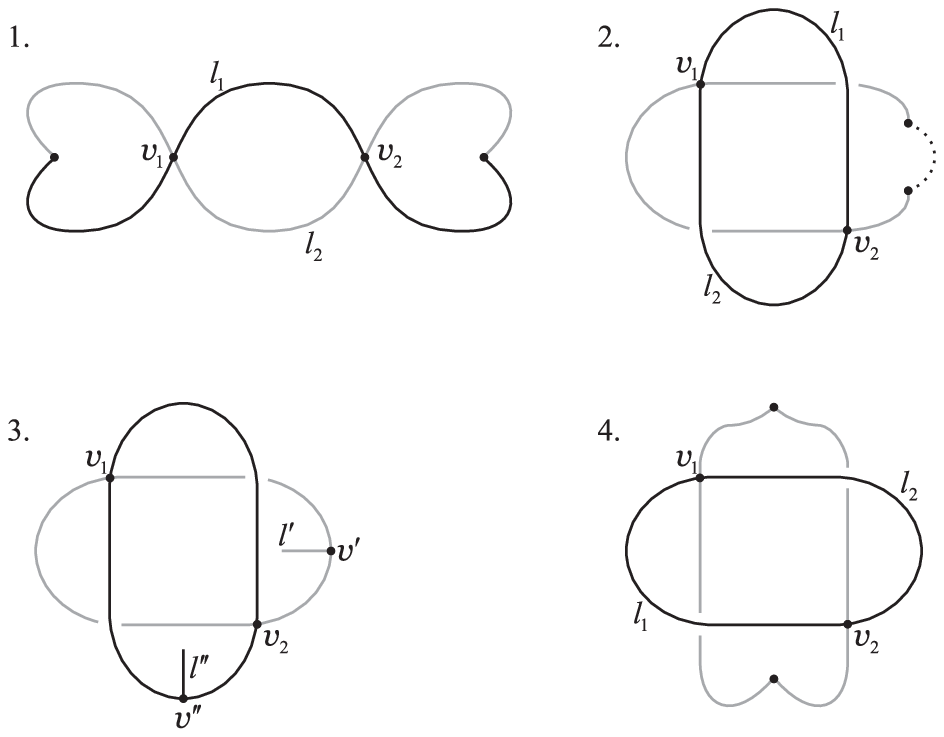}
    \caption{The four cases for $\widetilde G$ if $\lgh{ F}=6$.}
    \label{lgh_6_fold:fig}
    \end{center}
  \end{figure}

  Now, we consider the case where $\lgh{ F}=6$.
  As above, if $ F$ is incident to 3 different 4-valent vertices twice,
  we have that $S(\widetilde P)=\partial F$ (because
  $S(\widetilde P)$ is connected) and $\widetilde P$ has $3$ vertices,
  hence $\widetilde M$ is elliptic: a contradiction.
  So, if we consider the unfolded version of $ F$, we have a hexagon
  with four vertices identified in pairs (recall that $ F$ is incident
  to at most 4 distinct vertices).
  Using Lemmas~\ref{no_loop:lemma} and~\ref{no_self_adj_edge:lemma} as above, 
  we end up with 4 possible configurations,
  shown in Fig.~\ref{lgh_6_unfold:fig}.
  The corresponding portions of $\widetilde G$ adjacent to $ F$
  are shown in Fig.~\ref{lgh_6_fold:fig}.
  Each case is forbidden: let us show why.
  Cases~1,~2, and~4 lead to an embedded face $F'$ with $\lgh{F'}=2$, bounded by the loop
  $ l_1\cup l_2$ (where $
  F'\subset\widetilde\Sigma$ in case~1, and $
  F'\not\subset\widetilde\Sigma$ in cases~2 and~4), in contrast with Lemma~\ref{no_small_faces:lemma}.
  Concerning Case~3, the edges $l'$ and $ l''$ are different, 
  hence one of the two faces incident to $l'$ is
  incident to 5 different vertices (namely $ v_1$,
  $ v_2$, $ v'$, $ v''$, and the other
  endpoint of $ l'$), so we are done.

\begin{figure}[tp]
    \begin{center}
    \mettifig{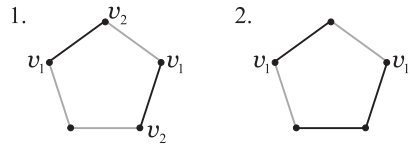}
    \caption{Two configurations when $\lgh{ F}=5$.
      In case~1 four vertices are identified in pairs, while in
      case~2 only two vertices are identified together.}
    \label{lgh_5_unfold:fig}
    \end{center}
  \end{figure}

\begin{figure}[tp]
    \begin{center}
    \mettifig{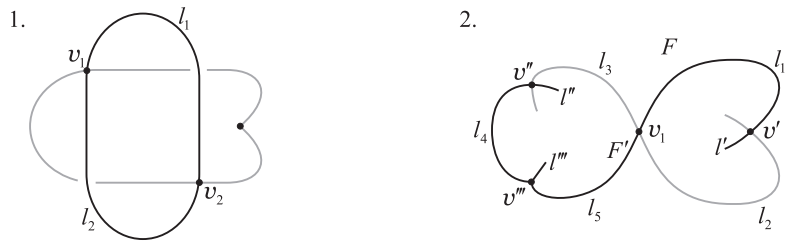}
    \caption{The two cases if $\lgh{ F}=5$.}
    \label{lgh_5_fold:fig}
    \end{center}
  \end{figure}

\begin{figure}[tp]
    \begin{center}
    \mettifig{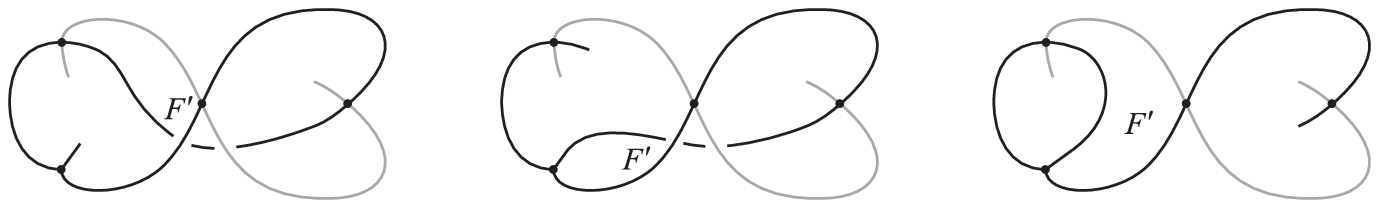}
    \caption{The three configurations of $\widetilde G$ if two of the
      three $ l^*$'s coincide, when $\lgh{ F}=5$.}
    \label{lgh_5_lgh_5:fig}
    \end{center}
  \end{figure}

\begin{figure}[tp]
    \begin{center}
    \mettifig{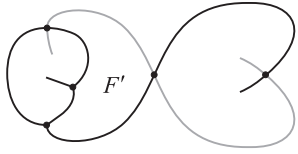}
    \caption{The (forbidden) configuration of $\widetilde G$ if $\lgh{
        F'}=4$.}
    \label{lgh_5_lgh_4:fig}
    \end{center}
  \end{figure}

  Finally, we consider the case where $\lgh{ F}=5$.
  The unfolded version of $ F$ is a pentagon
  with two or four vertices identified in pairs, 
  and using Lemmas~\ref{no_loop:lemma} and~\ref{no_self_adj_edge:lemma}
  we restrict ourselves to the two configurations
  drawn in Fig.~\ref{lgh_5_unfold:fig}, 
  yielding the two cases shown in Fig.~\ref{lgh_5_fold:fig}.
  Case~1 leads to an embedded face $F'$ (bounded by the loop
  $ l_1\cup l_2$ and non-contained in $\widetilde\Sigma$) with
  $\lgh{F'}=2$, in contrast with Lemma~\ref{no_small_faces:lemma}.
  Case~2 is slightly more complicated.
  Consider the face $F'$ shown in Fig.~\ref{lgh_5_fold:fig}-right.
  If two of the three $ l^*$'s coincide, we are done: in fact, either
  $F'$ is incident to $5$ distinct vertices, or $\lgh{ F'}>5$ (but
  $\lgh F$ is maximal), or $F'$ is an embedded triangle (contradicting
  Lemma~\ref{no_small_faces:lemma}), see Fig.~\ref{lgh_5_lgh_5:fig}.
  If instead the three $ l^*$'s are different, either $F'$ is incident
  to $5$ distinct vertices or it is a square as in
  Fig.~\ref{lgh_5_lgh_4:fig}.
  But such a square is excluded by Lemma~\ref{no_small_faces:lemma},
  and we are done.
\finedimo

\appendix
\section{On {\boldmath$I$}-bundles over tori and Klein bottles}\label{Ibundles:section}

Set $I = [-1,1]$.
In this appendix, we classify the $I$-bundles over the torus $T$ and the Klein bottle $K$. 
We denote by $D$, $A$, and $S$
respectively the disc, the annulus, and the M\"obius strip. 
Recall that there are two $S^1$-bundles $A\times S^1$ and $A\timtil S^1$
over $A$, and analogously two $S^1$-bundles $S\times S^1$ and $S\timtil S^1$
over $S$.
We denote by $\bar A$ the annulus with one mirror boundary
and by $\dot D$, $\ddot{D}$ respectively the disc with one or two mirror
segments in its boundary.
Therefore each $\bar A$ and $\dot D$ has one true boundary component,
while $\ddot{D}$ has two.
By mirroring one boundary component of $A\times S^1$ or $A\timtil S^1$
we get the two Seifert manifolds over the orbifold $\bar A$
(we denote them by $\bar A\times S^1$ and $\bar A\timtil S^1$).
Moreover, there is only one Seifert manifold over the orbifold
$\ddot{D}$: we denote it by $\ddot{D}\times S^1$.

\begin{prop} \label{I-bundles:prop}
There are, up to homeomorphism, two $I$-bundles $T\times I$ and $T\timtil I$ over $T$,
and three $I$-bundles $K\times I$, $K\timtil I$, $K\timtiltil I$ over $K$. 
We have $\partial (T\timtil I)\cong T$,
$\partial (K\timtil I)\cong T$, and $\partial (K\timtiltil I)\cong K$. 
They have the following Seifert fibrations:
\begin{itemize}
\item
  $T\times I$ fibers as $A\times S^1$,
\item
  $T\timtil I$ fibers as $S\times S^1$ and as $\bar A\times S^1$,
\item
  $K\times I$ fibers as $A\timtil S^1$ and as
  $\ddot{D}\times S^1$,
\item
  $K\timtil I$ fibers as $S\timtil S^1$ and as $\big(D,(2,1),(2,1)\big)$,
\item
  $K\timtiltil I$ fibers as $\bar A\timtil S^1$
  and as $\big(\dot D,(2,1)\big)$.
\end{itemize}
If $M$ is an $I$-bundle over $K$ or $T$ different from $K\timtil I$,
every fibration of one component of $\partial M$
extends to a Seifert fibration of $M$.
\end{prop}
\begin{proof}
The set of $I$-bundles over a closed surface $X$ up to fiber-preserving homeomorphisms 
is in 1-1 correspondence with
the orbits of $H^1(X;\matZ/_{2\matZ})$ under the action of the mapping class group
 of $X$. If $X=T$,
we have $H^1(T;\matZ/_{2\matZ}) = \matZ/_{2\matZ}\times \matZ/_{2\matZ}$, 
and using Dehn twists one sees that there are two orbits $\{(0,0)\}$ and 
$\{(1,0),(0,1),(1,1)\}$, giving respectively the product $T\times I$ 
and a non-orientable $I$-bundle, which we
denote by $T\timtil I$. If $X=K$, we have 
$H_1(K;\matZ) = \langle a,b|a+b=b+a,2a=0\rangle = \matZ/_{2\matZ}\times\matZ$
and $H^1(K;\matZ/_{2\matZ}) = {\rm Hom}(\matZ/_{2\matZ}\times\matZ,\matZ/_{2\matZ})=
\matZ/_{2\matZ}\times \matZ/_{2\matZ}$ again.
The mapping class group of $K$ is homeomorphic to 
$\matZ/_{2\matZ}\times \matZ/_{2\matZ}$ and is generated
by two automorphisms $\phi$ and $\psi$ whose action on $H_1(K,\matZ)$
is given by
$$\phi(a) =  a,\quad \phi(b)=-b \qquad {\rm and} \qquad \psi(a)=a, \quad \psi(b)=a+b.$$
(See the Appendix of~\cite{MaPe:nonori} for a proof of this fact.)
Therefore the orbits on $H^1(K;\matZ/_{2\matZ})$ are 
$\{(0,0)\}$, $\{(0,1)\}$, and $\{(1,0),(1,1)\}$, giving respectively
the (non-orientable) product $K\times I$, the orientable $K\timtil I$, and another non-orientable manifold
which we denote by $K\timtiltil I$.

Each such $I$-bundle can be described as a cube $I^3$ with the opposite lateral faces 
appropriately identified, so that
the horizontal $I^2\times\{0\}$ closes up to the zero-section.
The Seifert fibrations are then quotients of the two fibrations of $I^3$ 
given by $I\times\{p\}\times\{q\}$ and
$\{p\}\times I\times\{q\}$. We leave the details as an exercise for the reader.

Finally, let $N$ be an $I$-bundle, and let a component of $\partial N$ be fibered. 
If $N$ is a product, the fibration
extends trivially. If $N=K\timtiltil I$, the boundary $\partial N\cong K$ admits only two
non-isotopic fibrations, each of which extends to one of the two Seifert fibrations of $N$.
If $N=T\timtil I$, 
let $T_0$ be the zero-section. Fix generators $(\mu,\lambda)$ for $H_1(T_0;\matZ)$ so that
the $I$-bundle is determined by $\alpha\in H^1(T_0;\matZ/_{2\matZ})$
with $\alpha(\mu)=1$, $\alpha(\lambda)=0$.
Take also generators $(\mu',\lambda')$ for $H_1(\partial N;\matZ)$ which project to $(2\mu, \lambda)$.
With respect to these generators, every $\matr mnpq\in\GL_2(\matZ)$ with
even $n$ gives an automorphism of $T_0$ that
extends to an automorphism of $N$, acting as $\matr m{n/2}{2p}q$ on $\partial N$.
Via such automorphisms, an element of $H_1(\partial N;\matZ)$  is equivalent to either $\mu'$ or
$\lambda'$. Therefore, every given fibration of $\partial N$ is equivalent to
one of the two fibrations induced
by the two Seifert fibrations described above.
\end{proof}

The following corollary says that the two ``strange'' $I$-bundles
$T\timtil I$ and $K\timtiltil I$ do not occur near a Seifert block.
\begin{cor} \label{KtimtilI:cor}
  Let $M$ be a closed \ptwoirred\ manifold, and let $X$ be a $K$
  or a $T$ of the geometric decomposition of $M$.
  If $X$ is adjacent to a Seifert block,
  its regular neighborhood is either a product or $K\timtil I$.
\end{cor}
\begin{proof}
  The neighborhood is an $I$-bundle $N$ over $X$. 
  By Proposition~\ref{I-bundles:prop}, if $N$ is not a product or $K\timtil I$, then
  $\partial N$ is connected and the fibration of the adjacent Seifert block extends to $N$,
  a contradiction.
\end{proof}

\begin{cor} \label{negative:base:cor}
  Let $M$ be a closed \ptwoirred\ manifold with non-trivial geometric
  decomposition.
  Every Seifert block has hyperbolic base orbifold and finite volume.
\end{cor}
\begin{proof}
  Suppose a Seifert block has $\chiorb\geqslant 0$. If $\chiorb>0$, the base orbifold is either
  $D$ with one cone point at most, or $\dot D$. In those cases, the block is the solid torus or
  the solid Klein bottle, which is impossible.
  If instead $\chiorb=0$, the base orbifold
  is one of $A$, $\bar A$, $S$, $\ddot D$, $\dot D$ with one
  point with cone angle $\pi$, or $D$ with two points with cone angle
  $\pi$.
  There are two distinct Seifert fibrations over the orbifolds $A$,
  $\bar A$, and $S$, and one fibration
  over the other ones. By Proposition~\ref{I-bundles:prop}, the total
  space of every such fibration is homeomorphic
  to an I-bundle over $K$ or $T$. But no such block can occur in a
  geometric decomposition.
  Finally, note that the only blocks of the \JSJ\ decomposition with
  infinite volumes are flat.
\end{proof}

The following fact is not true for \JSJ\ decompositions.
\begin{cor} \label{decomposition:lifts:cor}
  Let $\widetilde M\to M$ be a finite covering of closed \ptwoirred\ manifolds.
  The pre-image of the geometric decomposition of $M$ is the geometric decomposition of $\widetilde M$.
\end{cor}
\begin{proof}
  A Seifert manifold with $\chiorb<0$ has a unique fibration~\cite{Sco}.
  Therefore, using Corollary~\ref{negative:base:cor}, we get that a 
  non-trivial decomposition is geometric if and only if
  every Seifert block has $\chiorb<0$ and the fibrations do not extend
  to any $K$ or $T$.
  Both conditions lift from $M$ to $\widetilde M$, hence we are done.
\end{proof}

\begin{cor} \label{even:cor}
Let $M$ be a non-orientable closed \ptwoirred\ manifold, whose geometric decomposition is made of Seifert blocks.
The geometric decomposition 
of $\widetilde M$ contains an even number of $K$'s.
\end{cor}
\begin{proof}
By Corollary~\ref{KtimtilI:cor}, the neighborhood of a $K$ in the geometric
decomposition of $M$ is either homeomorphic to $K\times I$ or 
to the orientable $K\timtil I$, giving rise respectively
to one $T$ or two $K$'s in the decomposition of $\widetilde M$.
\end{proof}

The following result is needed in the proof of Proposition~\ref{main:prop}.
\begin{cor} \label{Sol:cor}
A non-orientable manifold of \Sol\ geometry is a torus bundle over $S^1$, with some
monodromy $A\in\GL_2(\matZ)$ with $\det A=-1$. Two such manifolds with
monodromies $A, A'$ such that $|\tr A| = |\tr A'| \in\{1,2\}$ are homeomorphic.
\end{cor}
\begin{proof}
A manifold $M$ of \Sol\ geometry fibers over a $1$-orbifold, with
$T$'s and $K$'s as fibers.
If the $1$-orbifold is a segment (with two reflector endpoints), then
$M$ is the gluing of two $I$-bundles over $T$ or $K$ along
their connected boundaries.
Since $M$ is non-orientable, one $I$-bundle is either
$T\timtil I$ or $K\timtiltil I$, but in both cases $M$ is Seifert by Proposition~\ref{I-bundles:prop},
a contradiction. If instead the $1$-orbifold is $S^1$, then $M$ is a ($T$ or $K$)-bundle over $S^1$. But
$K$-bundles over $S^1$ are flat~\cite{Or}, hence it is a $T$-bundle, as required.

Suppose now we have two non-orientable manifolds with monodromies $A$ and $A'$ 
with $|\tr A| = |\tr A'| \in\{1,2\}$.
When $\det=-1$, we have $\matr abcd^{-1} = \matr {-d}bc{-a}$. Therefore we can suppose
$\tr A = \tr A' \in \{1,2\}$.
Taking $B\in\left\{\matr 10{\pm 1}1, \matr 1{\pm 1}01\right\}$ one sees
that $B^{-1}\matr abcd B \in \left\{\matr{a\pm b}b*{d\mp b}, \matr {a\mp c}*c{d\pm c}\right\}$.
Using this, we can suppose both $A$ and $A'$ have non-negative entries in the diagonal.
Therefore, if $\tr A = 1$ we get one $0$ entry in the diagonal and we easily get $A\sim A'$, whereas if
$\tr A = 2$ we either get one $0$ entry in the diagonal or $\matr 1121$, which is also easily transformed
into one having a $0$ entry in the diagonal. In both cases we get $A\sim A'$ and we are done.
\end{proof}

\section{On manifolds of type $\matH^2\times\matR$}

We prove here the following result.
\begin{prop} \label{small:Seifert:prop}
The two closed manifolds of type $\matH^2\times\matR$ with smallest base orbifold are
$$\big(\matRP^2,(2,1),(3,1)\big) \quad {\rm and} \quad \big(\bar D,(2,1),(3,1)\big).$$
Their base orbifold has volume $\pi/3$.
\end{prop}
\begin{proof}
The two manifolds above have $\chi^{\rm orb}=-\frac 16$, hence volume
$\frac {2\pi}6=\frac{\pi}3$.
We want to prove that every other manifold $M$ of type
$\matH^2\times\matR$ has $\chi^{\rm orb}<-\frac 16$, so let us suppose
by contradiction that $M$ is a manifold of type $\matH^2\times\matR$
with $\chi^{\rm orb}\geqslant -\frac 16$.

Let us first consider the case where the base orbifold of $M$ is $S^2$
with some $k$ points with cone angles $\frac{2\pi}{p_1},\ldots,$
$\frac{2\pi}{p_k}$.
Since $\chiorb=2-\sum(1-\frac 1{p_i})<0$, we have $k\geqslant 3$.
Suppose $k=3$. Then $\chi^{\rm orb} = 2-\sum(1-\frac 1{p_i}) = \sum \frac 1{p_i}-1$.
By our hypothesis $0>\sum\frac 1{p_i}-1\geqslant -\frac 16$, hence
$$(p_1,p_2,p_3)\;\;\in\;\;\{(2,3,h)\ |\ h\geqslant 7\}\cup\{(3,3,k)\ |\ 4\leqslant k\leqslant 
6\}\cup\{(2,4,l)\ |\ 5\leqslant l\leqslant 12\}.$$
We have that (the orientable) $M$ is $\big(S^2,(p_1,q_1),(p_2,q_2),(p_3,q_3)\big)$,
with Euler number $e=\sum_i\frac {q_i}{p_i}=0$. Therefore we get
$q_1p_2p_3+q_2p_3p_1+q_3p_1p_2 = 0$, hence $p_i|p_{i+1}p_{i+2}$
cyclically for $i=1,2,3$.
The only triple $(p_1,p_2,p_3)$ 
fulfilling this requirement is $(2,4,8)$.
But $q_3$ is odd, hence $4q_1+2q_2+q_3\neq 0$ gives $e\neq 0$ in that
case: a contradiction.

If $k=4$, we have $\chi^{\rm orb} = \sum\frac 1{p_i}-2$, 
hence $0>\sum\frac 1{p_i}-2\geqslant -\frac 16$. Therefore
$(p_1,p_2,p_3,p_4) = (2,2,2,3)$. Then $M=\big(S^2,(2,1),(2,1),(2,1),(3,q)\big)$
giving $e\neq 0$ again.
If $k\geqslant 5$, then $\chi^{\rm orb} \leqslant 2-\frac 52=-\frac
12$, and we are done.

If the orbifold is $\matRP^2$ or $\bar D$ with some $k$ points
with cone angles $\frac {2\pi}{p_1},\ldots,$ $\frac{2\pi}{p_k}$,
we have $\chi^{\rm orb} = 1 - \sum(1-\frac 1{p_i})$. 
Since $0>\chi^{\rm orb}\geqslant-\frac 16$, we get $k=2$ and
$(p_1,p_2)=(2,3)$, hence $M$ is one of the two listed above. 
Finally, if the surface underlying the base orbifold has $\chi \leqslant 0$, 
we get $\chi^{\rm orb} \leqslant -\frac 12$.
\end{proof}

\vspace{2cm}
\noindent
\hspace*{8cm}Dipartimento di Matematica\\ 
\hspace*{8cm}Universit\`a di Pisa\\ 
\hspace*{8cm}Via F. Buonarroti 2\\ 
\hspace*{8cm}56127 Pisa, Italy\\ 
\hspace*{8cm}amendola@mail.dm.unipi.it\\
\hspace*{8cm}martelli@mail.dm.unipi.it

\end{document}